\documentclass[11pt]{article}

\RequirePackage[OT1]{fontenc}

\RequirePackage[aap,amsthm,amsmath,natbib]{imsart}
\RequirePackage[aap]{imsart}
\usepackage{amsfonts}
\usepackage{amssymb}

\RequirePackage[dvips]{hyperref}

\startlocaldefs

\newcommand{\si}{\sigma}
\newcommand{\Si}{\Sigma}

\newcommand{\mbb}{\mathbb}
\newcommand{\dsty}{\displaystyle}
\newcommand{\dsum}{\sum\limits}
\newcommand{\dprod}{\prod\limits}
\newcommand{\dlim}{\lim\limits}

\newcommand{\ra}{\rightarrow}

\newcommand{\ol}{\overline}

\newcommand{\al}{\alpha}
\newcommand{\be}{\beta}
\newcommand{\de}{\delta}
\newcommand{\ep}{\epsilon}

\newcommand{\la}{\lambda}

\newcommand{\tri}{\triangle}
\newcommand{\itri}{\triangledown}
\theoremstyle{plain}
\newtheorem{thm}{Theorem}[section]

\theoremstyle{definition}

\theoremstyle{remark}

\endlocaldefs

\begin{document}

\begin{frontmatter}

\title{Large Deviations and Random Energy Models}
\runtitle{Random Energy Models}

\author{\fnms{N. K.} \snm{Jana}\ead[label=e2]{nabin\_r@isical.ac.in}}
\address{Department of Mathematics\\ Bijoy Krishna Girls' College\\ Howrah, India\\\printead{e2}}
\and
\author{\fnms{B. V.} \snm{Rao}\ead[label=e1]{bvrao@isical.ac.in}}
\address{Stat-Math Unit\\ Indian Statistical Institute\\ Kolkata, India \\\printead{e1}}
\affiliation{Indian Statistical Institute}

\runauthor{Jana and Rao}

\begin{abstract}
A unified treatment for the existence of free energy in several
random energy models is presented. If the sequence of distributions
associated with the particle systems obeys a large deviation
principle, then the free energy exists almost surely. This includes
all the known cases as well as some heavy-tailed distributions.
\end{abstract}



\end{frontmatter}
\section{Introduction}
The purpose of this note is to bring out the essence involved in the
existence theorems for the limiting free energy in several random
energy models. Basically, they are corollaries to the large
deviation principle (LDP) obeyed by certain empirical measures
coupled with Varadhan's lemma. Since this note is addressed to
probabilists, the results are first formulated in the large
deviation setting and then the consequences to the spin glass models
are explained in each section. At the first sight it may appear that
the large deviation principles proved here are nothing but the well
known principles for empirical measures, but however a little
reflection shows that it is not so. For large deviation terminology,
we refer to \citep{DE,DZ}.

We could have given the most general result of section 5 and then
deduced all the results from it. Instead, we decided to go from
simple model to the more general one, so that it is easier for the
reader to follow.

\section{Random Energy Model (REM)}

Let $\left(\la_N, N\geq1 \right)$ be a sequence of probabilities on
the real line $\mbb{R}$. Assume that $\la_N\Rightarrow \la$  and
satisfy large deviation principle with a convex rate function
$\mathcal{I}(x)$ on $\mbb{R}$. For every $N$, suppose $\xi_i, 1\leq
i\leq 2^N$ be i.i.d. random variables having distribution $\la_N$.
Define for each $\omega$, $\mu_N(\omega)$ to be the empirical
measure, namely $\mu_N(\omega)=\frac{1}{2^N}\sum
\delta(\xi_i(\omega))$.

\begin{thm}
For a.e. $\omega$ the above sequence $\{\mu_N(\omega)\}$ satisfies
LDP with rate function $\mathcal{J}$ given by,

$$\begin{array}{rll}
\mathcal{J}(x)& =\mathcal{I}(x) &\mbox{if }\quad \mathcal{I}(x)\leq \log2\\
&=\infty &\mbox{if }\quad \mathcal{I}(x)>\log2.
\end{array}$$

\end{thm}

\proof {\bf Step 1:} {\em Let $\tri$ be a subinterval of $\mbb{R}$.
If $\sum 2^N \la_N(\tri)<\infty$, then almost surely eventually
$\mu_N(\tri)=0.$}\\
Indeed, using $P$ for the probability on the space where the
random variables are defined,
$$P(\mu_N(\tri)>0)= P(\si_i \in \tri \mbox{ for some } i) \leq2^NP_N(\tri).$$
Use Borel - Cantelli.

{\bf Step 2:} {\em Let $\tri$ be a subinterval of $\mbb{R}$. If
$\sum \frac{1}{2^N \la_N(\tri)} <\infty$, then for any $\ep>0$
almost surely eventually $$(1-\ep) \la_N(\tri) \leq \mu_N(\tri) \leq
(1+\ep)\la_N(\tri).$$} Indeed, $$ \mbox{Var }\mu_N(\tri) =
E\left(\frac{1}{2^N}\sum 1_{\tri}(\xi)\right)^2 - \la_N^2(\tri) \leq
\frac{1}{2^N} \la_N(\tri).$$ Chebyshev yields
$$P(|\mu_N(\tri)-\la_N(\tri)|> \ep \la_N(\tri)) \leq \frac{1}{\ep^2
2^N \la_N(\tri)}.$$ Use Borel-Cantelli.\\

{\bf Step 3:} Fix an interval $\tri$ for which $\frac{1}{N}\log
\la_N(\tri)$ has a limit, say, $-L<-\log2$. Fix $\al>0$ such that
$-L<-\log2-\al$. Then for sufficiently large $N$, $\frac{1}{N}\log
\la_N(\tri) \leq -\log2-\al$, that is, $\la_N(\tri)\leq
2^{-N}e^{-N\al}$. In other words, $2^N\la_N(\tri) \leq e^{-N\al}$
for all large $N$. Thus by Step 1, a.s. eventually $\mu_N(\tri)=0$.
So $\frac{1}{N}\log \mu_N(\tri)=-\infty$ a. s. eventually.\\

{\bf Step 4:} Fix an interval $\tri$ for which $\frac{1}{N}\log
\la_N(\tri)$ has a limit, say, $-L>-\log2$. Fix $\al>0$ so that
$-L>-\log2 +\al$. So for large $N$, $\frac{1}{N}\log
\la_N(\tri)>-\log2 +\al$, that is, $\la_N(\tri)>2^{-N}e^{N\al}$. In
other words, $2^N\la_N(\tri)> e^{N\al}$. Now use Step 2 with
$\ep=\frac{1}{2}$ and take logarithms to get almost surely,
$$\lim_{N\ra\infty} \frac{1}{N}\log \mu_N(\tri) =-L =
\lim_{N\ra\infty} \frac{1}{N}\log \la_N(\tri)$$

Proof of the theorem is completed as follows. $\mathcal{I}$ being
convex, clearly the set $\{x: \mathcal{I}(x)=\log2\}$ has at most
two points. Let now $\tri$ be an open interval bounded away from
these two points. If $\tri \subset \{x: \mathcal{I}(x)<\log2\}$ then
$\lim_{N\ra\infty} \frac{1}{N} \log \mu_N(\tri) =-\inf_{x\in \tri}
\mathcal{I}(x)= -\mathcal{I}(\tri)$ by Step 4, where as, if
$\tri\subset \{x: \mathcal{I}(x)>\log2\}$ then $\lim_{N\ra\infty}
\frac{1}{N} \log \mu_N(\tri) = -\infty$ by Step 3. By theorem 4.1.11
in \citep{DZ}, $\mathcal{J}$ is the rate function
for the sequence $\{\mu_N(\omega)\}$ for all most every $\omega$.\\

The implications for REM\citep{D1} are clear. Here for fixed $N$,
one denotes the random variables $\xi$ by $H_N(\sigma)$ indexed by
$\sigma\in\{+1,-1\}^N$. One considers the partition function
$Z_N(\beta)=\dsum_{\sigma} e^{-\be H_N(\sigma)}=2^N E_\si e^{-\be
H_N}$ where $E_\si$ is the expectation with respect to the uniform
probability on the $\si$ space. The limit, $\dlim_{N\ra\infty}
\frac{1}{N} \log Z_N(\beta)$ when exists, is called the free energy
of the system.

Now suppose $\{\la_N\}$ is a sequence of probabilities as above and
$H_N(\sigma)=N\xi_N(\sigma)$ where for fixed $N$, $\xi_N(\sigma)$
are i.i.d. (as $\sigma$ varies) with distribution $\la_N$. Let us
assume that $\mathcal{J}$ has compact support. This is ensured by
assuming $\mathcal{I}$ to be a good rate function \citep{DZ} or at
least $\{x: \mathcal{I}(x)\leq \log2\}$ is bounded.

The case when $\la_N$ is centered Gaussian with  variance
$\frac{1}{N}$ translates to $H_N$ being of variance $N$. This is the
classical case \citep{T1,D1,DW}. Here
$\mathcal{I}(x)=\frac{x^2}{2}$. As a result, for almost every
$\omega$, the rate function of the sequence $\{\mu_N(\omega)\}$ is
$\mathcal{I}(x)=\frac{x^2}{2}$ if $\frac{x^2}{2} \leq \log 2$ and
$\infty$ otherwise. This in turn implies, by Varadhan's lemma
\citep{DE,DZ,JR1}, $\dlim_N \frac{1}{N} \log Z_N(\be)=\log 2\; -
\inf\limits_{x^2\leq 2\log 2}\{\beta x+ \frac{x^2}{2}\}$ which can
easily be evaluated.

The case when $\la_N$ has density $\frac{N}{2}e^{-N|y|}$ for
$-\infty<y<\infty$ corresponds to $H_N$ being two sided exponential
with parameter one considered in \citep{N1}. Weibull distributions
can also be considered \citep{E2,BBM,JR1}. It is clear that symmetry
of the distributions does not play any role. It is also easy to see
that when the random variables are, moreover, non negative then
there is no phase transition.

\section{Generalized Random Energy Model (GREM)}

The setup is the following. Let $n\geq 1$ be a fixed integer. For
each $j$, $1\leq j\leq n$, we have a sequence of probabilities
$\{\la_N^j, N\geq n\}$ on $\mbb{R}$ which weakly converges to
$\delta_0$ and obey LDP with a convex rate function $\mathcal{I}_j$.
Denote $\{+1,-1\}^N$ by $\Si_N$. For each $N$, let $k(1,N),\ldots,
k(n,N)$ be non-negative integers adding to $N$ and put
$\Si_{jN}=\{+1,-1\}^{k(j,N)}$. Clearly,
$\Si_N=\Si_{1N}\times\Si_{2N}\times\cdots\times\Si_{nN}$ and we
express $\sigma\in \Si_N$ as $\si_1\si_2\cdots\si_n$ with
$\si_i\in\Si_{iN}$, in an obvious way. For fixed $N$ we have a bunch
of independent random variables as follows: $\{\xi(\si_1):\,
\si_1\in\Si_{1N}\}$ having distributions $\la_N^1$,
$\{\xi(\si_1\si_2):\,\si_2\in\Si_{2N}, \si_1\in\Si_{1N}\}$ having
distributions $\la_N^2$ and in general
$\{\xi(\si_1\si_2\cdots\si_{j-1}\si_j):\, \si_j\in\Si_{jN},
\cdots,\si_1\in\Si_{1N}\}$ having distribution $\la_N^j$.

Define for each $\omega$, $\mu_N(\omega)$ to be the empirical
measure on $\mbb{R}^n$, namely,
$$\mu_N(\omega)=\frac{1}{2^N}\sum_\si
\delta\left<\xi(\si_1,\omega),\xi(\si_1\si_2,\omega),\cdots,\xi(\si_1\cdots\si_n,\omega)\right>$$
where $\delta\left<x\right>$ denotes the point mass at
$x\in\mbb{R}^n$.
\begin{thm}
Suppose $\frac{k(j,N)}{N}\ra p_j>0$ for $1\leq j\leq n$. Then for
a.e. $\omega$, the sequence $\{\mu_N(\omega), N\geq n\}$ satisfies
LDP with rate function $\mathcal{J}$ given as follows:

Supp$(\mathcal{J})=\{(x_1,\cdots,x_n): \dsum_{k=1}^j
\mathcal{I}_k(x_k) \leq \dsum_{k=1}^j p_k\log2\; \mbox{ for }\,
1\leq j\leq n \}$

and $$\begin{array}{rll} \mathcal{J}(x)& = \dsum_{k=1}^n
\mathcal{I}_k(x_k)& \mbox{if }
x\in \mbox{Supp}(\mathcal{J})\\
&=\infty&\mbox{otherwise}
\end{array}$$
\end{thm}
\proof The proof proceeds as in Theorem 2.1, we only explain the
steps involved. In what follows $\tri$ denotes a box in $\mbb{R}^n$
with sides $\tri_j;\; 1\leq j\leq n$ where each $\tri_j$ is an
interval.

{\bf Step 1:} {\em If for some $j$, $$\sum_{N\geq
n}\prod_{i=1}^j2^{k(i,N)}\la_N^i(\tri_i)<\infty$$ then a.s.
eventually $\mu_N(\tri)=0$.}\\
{\bf Step 2:} {\em If for each $j$, $$\sum_{N\geq n}
\prod_{i=1}^j\frac{1}{2^{k(i,N)}\la_N^i(\tri_i)}<\infty$$ then for
any $\ep>0$ a.s. eventually
$$(1-\ep)\prod_{i=1}^n\la_N^i(\tri_i) \leq \mu_N(\tri)\leq
(1+\ep)\prod_{i=1}^n\la_N^i(\tri_i).$$}

This step involves calculation of $\mbox{Var }(\mu_N(\tri))$ which
is carried out in a more general set up in section 5. The reader may
also consult \citep{CCP,DD,JR1}. The remaining two steps are
accordingly modified.

The implications for GREM\citep{D2} are clear. For fixed $N$, and
$\si\in\Si_N$ one defines the Hamltonian
$$H_N(\si)=N\sum_{i=1}^n a_i \xi(\si_1\cdots\si_i).$$ Here $a_i,\;
1\leq i\leq n$ are positive numbers called weights. In the Gaussian
case, it is customary to take $\sum a_i^2=1$, though it is not a
mathematical necessity. As earlier, $Z_N(\be)=\sum_\si e^{-\be
H_N(\si)}$. Special choices of $\la_N^i$ lead to all the known
models considered. Centered Gaussians were consider in
\citep{D2,CCP,DD,JR1}. More general distributions as well as the
cases when some $p_j$ are zero were considered in \citep{JR1}.
Moreover one could take different distributions for different values
of $j$, see \citep{JR2} for some interesting consequences. Thus the
main problem of GREM is reduced to a variational problem. Note that,
if $n=1$, GREM reduces to REM.

\section{Bolthausen-Kistler Model (BKM)}

We follow the same notation as in the previous section for $\Si_N,
\Si_{iN}, k(i,N)$ and $\si=\si_1\si_2\cdots\si_n$. Let $I=\{1,
2,\cdots,n\}$ and $S$ be the collection of non-empty subsets of $I$.
Suppose that for each non-empty $s\subseteq I$, we have a sequence
of probabilities $\{\la_N^s:\; N\geq n\}$ weakly converging to
$\de_0$ and obeying LDP with a convex rate function
$\mathcal{I}_s(x)$. For $s=\{ i_1, i_2, \cdots,i_k\}\in S$ where
$i_1< i_2< \cdots<i_k$, we denote $\si(s)=\si_{i_1} \si_{i_2}\cdots
\si_{i_k}$. Now for fixed $N$, we have a bunch of independent random
variables $\xi(s,\si(s))$ as $s$ and $\si(s)$ vary. For $s\in S$,
all the $\xi(s,\si(s))$ have distribution $\la_N^s$. We define the
empirical measure on $\mbb{R}^{S}$,  by
$$\mu_N(\omega)=\frac{1}{2^N}\sum_\si \delta\left<\xi(s,\si(s),\omega):\;
s\in S\right>.$$

\begin{thm}
Suppose $\frac{k(j,N)}{N}\ra p_j>0$ for $1\leq j\leq n$. Then for
a.e. $\omega$, the sequence $\{\mu_N(\omega), N\geq n\}$ satisfies
LDP with rate function $\mathcal{J}$ given as follows:

Supp$(\mathcal{J})=\{(x_s,s\in S):\; \forall s\in S,\;
\dsum_{t\subseteq s} \mathcal{I}_t(x_t) \leq \dsum_{k\in s} p_k\log2
\}$

and $$\begin{array}{rll} \mathcal{J}(x)& = \dsum_{s\in S} \mathcal{I}_s(x_s)& \mbox{if }
x=(x_s)\in \mbox{Supp}(\mathcal{J})\\
&=\infty&\mbox{otherwise}
\end{array}$$
\end{thm}

The proof proceeds along the same lines as in Theorem 3.1. The only
complications, again, are in calculating $\mbox{Var } \mu_N(\tri)$.
See next section.

In BKM\citep{BK}, we have Hamiltonian $H_N(\si)=N\dsum_{s\in S} a_s
\xi(s, \si(s))$ where $a_s, \; s\in S$ are non-negative weights.
Again if $\mathcal{J}$ has compact support, Varadhan's lemma reduces
calculation of free energy to that of a variational problem. In BKM,
they do not consider all non-empty subsets of $I$, rather, some sub
collection. But that can be achieved by taking some weights $a_s$ to
be zero appropriately. Further, if one considers a chain like
$\{1\}, \{1,2\},\cdots,\{1,2,\cdots,n\}=I$ and takes non zero
weights $a_s$ for $s$ only in the chain, one gets GREM.

\section{Models with External Field}

The model described in this section includes all the above and also
incorporates external field. To describe the set up, we will use the
same notation as above except for $S$ and $s$. Here the set of all
non-empty ordered sequences $s=\left< i_1, i_2,\cdots, i_k\right>$
of distinct elements from $I$, having length $\leq n$ will be
denoted as $S$. We use the same notation of $\si(s)$ as in the
previous section. For $s\in S$, we have a sequence of probabilities
$\{\la_N^s,\, N\geq n\}$ weakly converging to $\de_0$ obeying LDP
with a convex rate function $\mathcal{I}_s$. As earlier, for fixed
$N$, we have a bunch of independent random variables $\xi(s,\si(s))$
as $s$ and $\si(s)$ vary. For $s\in S$, all the $\xi(s,\si(s))$ have
distribution $\la_N^s$. Let $\overline{\si}_i$ denote the sum of the
$k(i,N)$ many $+1$ and $-1$ appearing in $\si_i$. We define the
empirical measure on $\mbb{R}^{S}\times\mbb{R}^n$ by
$$\mu_N(\omega)=\frac{1}{2^N}\sum_\si \delta\left<\left(\xi(s,\si(s),\omega):\;
s\in S\right),\left(\frac{\overline{\si}_i}{N}:\;1\leq i \leq
n\right)\right>.$$ Points in $\mbb{R}^{S}\times\mbb{R}^n$ will be
denoted by $\left((x_s, s\in S),(y_i, i\leq n)\right)$ or simply as
$(x_{_S},y_{_I})$. For $A\subset I$, we will denote $S_A$  to be the
set of non-empty ordered sequences of distinct elements from $A$.

First note that, by Cramer's theorem \citep{DZ}, the arithmetic
averages of i.i.d. mean zero, $\pm 1$ valued random variables
satisfy LDP with rate function $\mathcal{I}_0$ where
$\mathcal{I}_0(y)=\infty$ for $|y|>1$ and for $-1\leq y\leq 1$,
$$\begin{array}{rl}\mathcal{I}_0(y)&= y\tanh^{-1}y
-\log\cosh(\tanh^{-1}y)\vspace{1ex}\\
&=\frac{1+y}{2}\log(1+y)+\frac{1-y}{2}\log(1-y). \end{array}$$

\begin{thm}
Suppose $\frac{k(j,N)}{N}\ra p_j>0$ for $1\leq j\leq n$. Then for
a.e. $\omega$, the sequence $\{\mu_N(\omega), N\geq n\}$ satisfies
LDP with rate function $\mathcal{J}$ given as follows:

Supp$(\mathcal{J})=\left\{(x_{_S},y_{_I}):\; \forall A\subset I,
\dsum_{t\in S_A} \mathcal{I}_t(x_t)+\dsum_{i \in A}
p_i\mathcal{I}_0\left(\frac{y_i}{p_i}\right) \leq \dsum_{k\in A}
p_k\log2 \right\}$

and $$\begin{array}{rll} \mathcal{J}(x_{_S},y_{_I})& = \dsum_{s\in
S} \mathcal{I}_s(x_s) + \dsum_{i \in I}
p_i\mathcal{I}_0\left(\frac{y_i}{p_i}\right)& \mbox{if }
(x_{_S},y_{_I})\in \mbox{Supp}(\mathcal{J})\\
&=\infty&\mbox{otherwise}
\end{array}$$
\end{thm}

\proof In what follows, $\square=\prod_{s\in
S}\tri_s\times\prod_{i=1}^n\itri_i$ is a box in
$\mbb{R}^S\times\mbb{R}^n$ where $\tri_s$ for each $s\in S$ and
$\itri_i$ for $i\leq n$ are subintervals of $\mbb{R}$. Also $A$
denotes non empty subset of $I$. For $A\subseteq I$, we denote
$Q_{AN}=\dprod_{s\in S_A} \la_N^s(\tri_s)$; $k(A,N)=\dsum_{i\in A}
k(i,N)$ and $\al_{AN}=\frac{1}{2^{k(A,N)}}\dsum_{\left<\si_i:\; i\in
A\right>}\dprod_{i\in
 A}1_{\itri_i}(\frac{\overline{\si}_i}{N})$.

{\bf Step 1:} {\em If for some $A\subseteq I$, $\sum_{N\geq n}
2^{k(A,N)} Q_{AN} \al_{AN}<\infty$ then almost surely eventually
$\mu_N(\square)=0$.} \\

Let $A$ be such that $\sum_{N\geq
n}2^{k(A,N)}Q_{AN}\al_{AN}<\infty$. Then
$$\begin{array}{rl}
\mu_N(\square)&=\frac{1}{2^N}\dsum_\si \dprod_{s\in S}
1_{\tri_s}\left(\xi(s,\si(s))\right)\dprod_{i\leq
n}1_{\itri_i}\left(\frac{\ol{\si}^i}{N}\right)\vspace{1ex}\\

&\leq\frac{1}{2^N}\dsum_\si \dprod_{s\in \mathcal{S}_A}
1_{\tri_s}\left(\xi(s,\si(s))\right)\dprod_{i\in A}
1_{\itri_i}\left(\frac{\ol{\si}^i}{N}\right)\vspace{1ex}\\

&=\frac{1}{2^{k(A,N)}}\dsum_{\left<\si_i : i\in A\right>}
\dprod_{s\in \mathcal{S}_A}
1_{\tri_s}\left(\xi(s,\si(s))\right)\dprod_{i\in A}
1_{\itri_i}\left(\frac{\ol{\si}^i}{N}\right)\\
\end{array}$$

As a consequence,
$$\begin{array}{rl}
P(\mu_N(\square)>0)&=P\left(\dsum_{\left<\si_i : i\in A\right>}
\dprod_{s\in \mathcal{S}_A}
1_{\tri_s}\left(\xi(s,\si(s))\right)\dprod_{i\in
A}1_{\itri_i}\left(\frac{\ol{\si}^i}{N}\right)\geq 1\right)\vspace{1ex}\\

&\leq Q_{AN}\dsum_{\left<\si_i : i\in A\right>} \dprod_{i\in
A}1_{\itri_i}\left(\frac{\ol{\si}^i}{N}\right)\vspace{1ex}\\

&=2^{k(A,N)}Q_{AN}\al_{AN}.
\end{array}$$

The hypothesis and Borel-Cantelli lemma complete the proof of the
claim.\\

{\bf Step 2:} {\em If for all $A\subseteq I$, $\sum_{N\geq
n}2^{-k(A,N)}Q_{AN}^{-1}\al_{AN}^{-1} < \infty$ then for any
$\ep>o$, almost surely eventually, $$(1-\ep)Q_{IN}\al_{IN}\leq
\mu_N(\square)\leq(1+\ep)Q_{IN}\al_{IN}.$$}

Note that
$$\begin{array}{l}
\mbox{Var}(\mu_N(\square))\\

=\frac{1}{2^{2N}}\dsum_\si\dsum_\tau E\left( \dsty\prod_{s\in S}
1_{\tri_s}\left(\xi(s,\si(s))\right)1_{\tri_s}\left(\xi(s,\tau(s))\right)\right)\dsty\prod_{i\leq
n}1_{\itri_i}\left(\frac{\ol{\si}^i}{N}\right)1_{\itri_i}\left(\frac{\ol{\tau}^i}{N}\right)\vspace{1ex}\\

\hspace{60ex}-Q_{IN}^2\al_{IN}^2\vspace{2ex}\\

\leq\frac{1}{2^{2N}}\dsum_{A\subseteq
I}\dsum_\si\dsum_{\tiny\begin{array}{c}\tau_i=\si_i, \forall i\in
A\\\tau_i\neq \si_i, \forall i\in A^c
\end{array}}E\left( \dsty\prod_{s\in S}
1_{\tri_s}\left(\xi(s,\si(s))\right)1_{\tri_s}\left(\xi(s,\tau(s))\right)\right)\times\vspace{1ex}\\

\hspace{50ex}\dsty\prod_{i\leq
n}1_{\itri_i}\left(\frac{\ol{\si}^i}{N}\right)1_{\itri_i}\left(\frac{\ol{\tau}^i}{N}\right)\vspace{2ex}\\

=\frac{1}{2^{2N}}\dsum_{A\subseteq
I}\frac{Q_{IN}^2}{Q_{AN}}\dsum_\si\dsum_{\tiny\begin{array}{c}\tau_i=\si_i,
\forall i\in A\\\tau_i\neq \si_i, \forall i\in A^c
\end{array}}\dsty\prod_{i\in A}1_{\itri_i}\left(\frac{\ol{\si}^i}{N}\right)\dsty\prod_{i\in
A^c}1_{\itri_i}\left(\frac{\ol{\si}^i}{N}\right)1_{\itri_i}\left(\frac{\ol{\tau}^i}{N}\right)\vspace{2ex}\\

=\dsum_{A\subseteq I}
\frac{Q_{IN}^2}{Q_{AN}}\,\frac{1}{2^{k(A,N)}}\,\frac{\al_{IN}^2}{\al_{AN}}.
\end{array}$$

Here is how we get the inequality above. The terms appearing in the
first expression corresponding to $\si_i \neq \tau_i$ for all $i\in
I$ are canceled with $Q_{IN}^2\al_{IN}^2$. Now by Chebyshev's
inequality
\[P(|\mu_N(\square)-E\mu_N(\square)|>\ep E\mu_N(\square))<
\frac{1}{\ep^2}\dsum_{A\subseteq I}
\frac{1}{2^{k(A,N)}Q_{AN}\al_{AN}}.\] Once again Borel-Cantelli
lemma and the hypothesis yield that a.s. eventually,
$$(1-\ep)E\mu_N(\square) \leq \mu_N(\square) \leq (1+\ep)
E\mu_N(\square).$$
Since $E\mu_N(\square)=Q_{IN}\al_{IN}$, proof of the claim is complete. \\

{\bf Step 3:} If a box $\square\subset\mbb{R}^{S}\times\mbb{R}^n$
has empty intersection with $Supp(\mathcal{J})$, then by Step 1,
almost surely eventually $\mu_N(\square)=0$ and hence
$\lim\limits_{N\ra \infty}\frac{1}{N} \log \mu_N(\square) =
-\infty.$

{\bf Step 4:} Suppose that
$\square\subset\mbb{R}^{S}\times\mbb{R}^n$ has non-empty
intersection with the interior of $Supp(\mathcal{J})$. Then by Step
2, almost surely, $\lim\limits_{N\ra \infty}\frac{1}{N} \log
\mu_N(\square) = \lim\limits_{N\ra \infty}\frac{1}{N} \log Q_{IN}
\al_{IN} = -\left[\dsum_{s\in S} \mathcal{I}_s(\tri_s)+\dsum_{1\leq
i \leq n}p_i\mathcal{I}_0(\frac{1}{p_i}\itri_i)\right]$. To see
this, one has only to decipher $Q_{IN}\al_{IN}$.

Since there are enough boxes $\square$ -- by convexity of
$\mathcal{I}_s$ -- satisfying conditions of either Step 3 or Step 4
to form a basis, the proof is completed as in Theorem 2.1 using
Theorem 4.1.11 in \citep{DZ}.

Implication of the above result will be clear if one defines the
Hamiltonian $H_N(\si)=N\dsum_{s\in S} a_s \xi(s, \si(s)) + h
\ol{\si}$ where $a_s, \; s\in S$ are non-negative weights and $h>0$
is the strength of the external field. Here also if $\mathcal{J}$
has compact support (which is true when the rate functions
$\mathcal{I}_s$ are convex), Varadhan's lemma reduces calculation of
free energy to that of a variational problem. Of course, it is not
always possible to solve this variational problem to arrive at a
closed form expression. When the external field is 0 and when only
increasing sequences have positive weight, then this reduces to BKM.

\section{Remarks}

1. As the reader would have noticed, convexity of the rate function
$\mathcal{I}$ is not essential. For example, the results are valid
if the rate function strictly decreases on $(-\infty, 0)$ and
strictly increasing on $(0, \infty)$. In fact, this will then
include the Weibull distribution with shape parameter smaller than
one.

2. The factor $2^N$ could easily be replaced by $k^N$ ($k$ integer
$\geq1$) with appropriate changes in the theorems.

3. We do not have any precise formulation of the results for GREM
with $n=\infty$ \citep{CCP}. Perhaps these arguments do not work for
SK-model.

\end{document}